\documentclass{elsart}
\usepackage{amssymb,amsfonts}

\newcommand{\trace}{\mathop{\rm Tr}\nolimits}

\newcommand{\diag}{\mathop{\rm Diag}\nolimits}

\newcommand{\twomat}[4]{\left(\begin{array}{cc}#1&#2\\#3&#4\end{array}\right)}

\newcommand{\cP}{{\mathcal P}}

\newcommand{\R}{{\mathbb{R}}}

\DeclareRobustCommand\openone{\leavevmode\hbox{\small1\normalsize\kern-.33em1}}
\newcommand{\id}{\mathrm{\openone}}

\newcommand{\ddt}{\left.\frac{\partial}{\partial t}\right|_{t=0}}

\newcommand{\be}{\begin{equation}}
\newcommand{\ee}{\end{equation}}
\newcommand{\bea}{\begin{eqnarray}}
\newcommand{\eea}{\end{eqnarray}}
\newcommand{\beas}{\begin{eqnarray*}}
\newcommand{\eeas}{\end{eqnarray*}}

\newtheorem{lemma}{Lemma}
\newtheorem{corollary}{Corollary}

\newtheorem{proposition}{Proposition}
\newtheorem{question}{Question}

\newcount\minute
\newcount\hour
\def\currenttime{%
    \minute\time
    \hour\minute
    \divide\hour60
    \the\hour:\multiply\hour60\advance\minute-\hour\the\minute}
\begin{document}
\begin{frontmatter}
\title{On Ando's inequalities for convex and concave functions}
\author{Koenraad M.R.\ Audenaert}
\address{
Institute for Mathematical Sciences, Imperial College London,\\
53 Prince's Gate, London SW7 2PG, United Kingdom}
\ead{k.audenaert@imperial.ac.uk}
\author{Jaspal Singh Aujla}
\address{Department of Applied Mathematics, National Institute of Technology,\\
Jalandhar 144011, Punjab, India}
\ead{aujlajs@nitj.ac.in}
\date{\today, \currenttime}
\begin{abstract}
For positive semidefinite matrices $A$ and $B$,
Ando and Zhan proved the inequalities
$||| f(A)+f(B) ||| \ge ||| f(A+B) |||$ and
$||| g(A)+g(B) ||| \le ||| g(A+B) |||$, for any unitarily invariant norm, and
for any non-negative operator monotone $f$ on $[0,\infty)$ with inverse function $g$.
These inequalities have very recently been generalised to non-negative concave functions
$f$ and non-negative convex functions $g$, by Bourin and Uchiyama, and Kosem, respectively.

In this paper we consider the related question whether the inequalities
$||| f(A)-f(B) ||| \le ||| f(|A-B|) |||$,
and
$||| g(A)-g(B) ||| \ge ||| g(|A-B|) |||$,
obtained by Ando, for operator monotone $f$ with inverse $g$, also
have a similar generalisation to non-negative concave $f$ and convex $g$.
We answer exactly this question, in the negative for general matrices, and affirmatively
in the special case when $A\ge ||B||$.

In the course of this work, we introduce the novel notion of $Y$-dominated majorisation between the spectra
of two Hermitian matrices, where $Y$ is itself a Hermitian matrix, and prove a certain property of this
relation that allows to strengthen the results of Bourin-Uchiyama and Kosem, mentioned above.
\end{abstract}

\end{frontmatter}
\section{Introduction}
In \cite{ando88}, Ando proved the following inequalities for
positive semidefinite (PSD) matrices $A$, $B$, and any unitarily
invariant (UI) norm.
For any non-negative operator monotone function $f(t)$ on $[0,\infty)$:
\be\label{eq:A1}
||| f(A)-f(B) ||| \le ||| f(|A-B|) |||,
\ee
and, when $f(0)=0$ and $f(\infty)=\infty$, and $g$ is the inverse function of $f$,
\be\label{eq:A2}
||| g(A)-g(B) ||| \ge ||| g(|A-B|) |||.
\ee
In a later paper \cite{andozhan}, Ando and Zhan proved the related
inequalities (with the same conditions on $f$ and $g$)
\be\label{eq:AZ1}
||| f(A)+f(B) ||| \ge ||| f(A+B) |||,
\ee
\be\label{eq:AZ2}
||| g(A)+g(B) ||| \le ||| g(A+B) |||.
\ee
The conditions on $f$ are satisfied by every operator concave
function $f$ with $f(0)=0$.

Inequality (\ref{eq:AZ2}) was generalised by Kosem \cite{kosem}
to non-negative convex functions $g$ on $[0,\infty)$, with $g(0)=0$.
Inequality (\ref{eq:AZ1}) was generalised very recently to any
non-negative concave function on $[0,\infty)$ by Bourin and
Uchiyama \cite{bourin}, who also gave a simpler proof of
Kosem's result.

In the same spirit, we consider the question whether inequalities
(\ref{eq:A1}) and (\ref{eq:A2}) can also be generalised to non-negative
concave $f$ and convex $g$, respectively.
After introducing the necessary prerequisites in Section \ref{sec:pre}, we present our main results
concerning this question in Section \ref{sec:main}. Regrettably, most of our results are negative answers,
and we give counterexamples to this generalisation.
The answer is even negative for the special case $A\ge B$, although the apparent hardness of finding counterexamples
had led us temporarily into believing that in that case the generalisation might actually hold.

All is not bad news, however. In Section \ref{sec:evidence} we answer the question
affirmatively in the special case when $A\ge ||B||$.
In Section \ref{sec:domi}, we introduce the novel notion of $Y$-dominated majorisation between the spectra
of two Hermitian matrices, where $Y$ is itself a Hermitian matrix. We prove a certain property of this
relation, namely Proposition \ref{prop:4b}, which we subsequently use, first in a rather destructive fashion.
To wit, the Proposition has been instrumental in finally discovering a counterexample to
the generalisation of (\ref{eq:A1}) for $A\ge B$; this will be reported in Section \ref{sec:counter}.
On the more constructive side, the Proposition also allows to strengthen the results of Bourin-Uchiyama and Kosem
mentioned above. This is the topic of the final Section, along with a few other applications.
\section{Preliminaries\label{sec:pre}}
In this Section, we introduce the notations and necessary prerequisites;
a more detailed exposition can be found, e.g.\ in \cite{bhatia}.
We will use the abbreviations LHS and RHS for left-hand side and right-hand side, respectively.

We denote the vector of diagonal entries of a matrix $A$ by $\diag(A)$.

We denote the absolute value by $|\cdot|$, both for scalars and for matrices. For matrices this is defined
as $|A| := (A^* A)^{1/2}$. Similarly, we denote the positive part of a real scalar or Hermitian
matrix by $(\cdot)^+$, and define it by $A^+ := (A+|A|)/2$.

In this paper, we are mainly concerned with monotonously increasing convex and concave functions from $\R$ to $\R$.
Kosem noted in \cite{kosem} that any such function can be approximated by a sum of angle functions
$x\mapsto ax+b(x-x_0)^+$, where $a\ge0$, and $b>0$ for a convex angle function ($b<0$ for a concave one).

We are also concerned with the unitarily invariant (UI) matrix norms, which we denote by $|||\cdot|||$, and which
are defined in terms of the singular values $\sigma_j(\cdot)$ of the matrix only.
We adopt the customary convention that the singular values are sorted in non-increasing order:
$\sigma_1\ge \sigma_2\ge\ldots\ge\sigma_d$.
Special cases of these norms are the operator norm $||\cdot||$, which is just equal to the
largest singular value $\sigma_1(\cdot)$,
and the Ky Fan norms $||\cdot||_{(k)}$, which are defined as the sum of the $k$ largest singular values:
$$
||A||_{(k)} := \sum_{j=1}^k \sigma_j(A).
$$

The famous Ky Fan dominance theorem states that a matrix $B$ dominates another matrix $A$ in all UI norms
if and only if it does so in all Ky Fan norms.
The latter set of relations can be written as a weak majorisation relation between the vectors of
singular values of $A$ and $B$:
$$
\sigma(A)\prec_w \sigma(B): \qquad \sum_{j=1}^k \sigma_j(A)\le \sum_{j=1}^k \sigma_j(B), \forall k.
$$
For PSD matrices, the above domination relation translates to a weak majorisation between
the vectors of eigenvalues: $\lambda(A)\prec_w\lambda(B)$.

Weyl's monotonicity Theorem (\cite{bhatia}, Corollary III.2.3) states that
$$
\lambda^\downarrow_k(A) \le \lambda^\downarrow_k(A+B), \quad\forall k,
$$
for Hermitian $A$ and positive semi-definite $B$.
Here, $\lambda^\downarrow(A)$ denotes the (real) vector of eigenvalues of $A$ sorted in non-increasing order.

Finally, we refer the reader to Chapter 2 of \cite{kato} for an exposition of a number of important functional analytic
properties of eigenvalues and corresponding eigenspaces of a Hermitian matrix, which we will need
in the proof of Proposition \ref{prop:delta}.
\section{Main Results\label{sec:main}}
The question we start with is about the straightforward generalisation of inequality (\ref{eq:A2}) to
non-negative convex functions.
\begin{question}\label{th2a}
For all $A,B,\ge0$, for all UI norms, and for non-negative convex functions
$g$ on $[0,\infty)$ with $g(0)=0$, does the inequality
$||| g(A)-g(B) ||| \ge ||| g(|A-B|) |||$ hold?
\end{question}
The answer to this question is negative, as shown by the following counterexample.
We consider the convex angle function $g(x) = x+(x-1)^+$ and the operator norm.
For the $2\times 2$ PSD matrices
\be
A = \left(\begin{array}{rr}
0.9 &0 \\
0& 0.6
\end{array}
\right),\quad
B = \left(
\begin{array}{rr}
0.8 &0.5 \\
0.5& 0.4
\end{array}
\right),
\ee
the eigenvalues of $g(|A-B|)$ are $0.65249$ and $0.35249$,
while those of $g(A)-g(B)$ are $0.65010$ and $-0.48862$. Thus,
$||g(|A-B|)||_\infty =0.65249$, which is larger than
$||g(A)-g(B)||_\infty=0.65010$.
\qed

Under the additional restriction $A\ge B$, the absolute value in the argument
of $g$ in the right-hand side vanishes, leading to a simplified statement, and
a second question, with better hopes for success. Introducing the matrix $\Delta=A-B$,
\begin{question}\label{th2}
For all $B,\Delta\ge0$, for all UI norms, and for non-negative convex functions
$g$ on $[0,\infty)$ with $g(0)=0$, does the inequality
$||| g(B+\Delta)-g(B) ||| \ge ||| g(\Delta) |||$ hold?
\end{question}
This restricted case also turns out to have a negative answer.
Counterexamples, however, were much harder to find, and required a reduction of the problem
based on certain results about a novel majorisation-like relation, which we call the
$Y$-dominated majorisation.
This will be the subject of Sections \ref{sec:domi} and \ref{sec:counter},
where a number of Propositions of independent interest are proven.

\bigskip

It is also very reasonable to ask:
\begin{question}
For all $B,\Delta\ge0$, for all UI norms, and for non-negative concave functions
$f$ on $[0,\infty)$, does the inequality
$||| f(B+\Delta)-f(B) ||| \le ||| f(\Delta) |||$ hold?
\end{question}
In fact, if this were true, a positive answer to Question \ref{th2} would easily follow,
using the same reasoning that was used in \cite{bourin} to derive the generalisation
of (\ref{eq:AZ1}) from the generalisation of (\ref{eq:AZ2}).

Again, this statement is false, as the following counterexample shows.
Consider the concave angle function $f(x) = \min(x,1) = x-(x-1)^+$,
and the $3\times 3$ PSD matrices
$$
B = \left(
\begin{array}{rrr}
         0.701816  &       0.317887   &      0.198910\\
         0.317887  &       1.014950   &     -0.093826\\
         0.198910  &      -0.093826   &      0.274236
\end{array}
\right)
$$
and
$$
\Delta = \left(
\begin{array}{rrr}
         0.192713 & 0   &      0 \\
         0  &       0.446505 & 0 \\
         0  &       0   &      0.455416
\end{array}
\right).
$$
One gets
$$
||f(\Delta)||_\infty = 0.455416
$$
while
$$
||f(B+\Delta)-f(B)||_\infty = 0.455776.
$$
\qed

\bigskip

In Section \ref{sec:evidence}, we consider an even more restricted special case, in which
the inequalities (\ref{eq:A1}) and (\ref{eq:A2}) finally do hold.
We actually prove that a stronger relationship holds in this special case:
\begin{proposition}\label{prop:ggc}
For non-negative, monotonously increasing and concave functions $g$, and $A,B\ge 0$ such that $A\ge ||B||$,
we have
\be\label{eq:ggc}
\lambda^\downarrow(g(A-B)) \ge \lambda^\downarrow(g(A)-g(B)).
\ee
\end{proposition}
An easy Corollary is the corresponding statement for monotonously increasing convex functions.
\begin{corollary}\label{prop:gg}
Let $f$ be a non-negative convex function on $[0,\infty )$ with $f(0)=0$.
Let $A,B\ge 0$ such that $A\ge ||B||$.
Then
\be\label{eq:gg}
\lambda^\downarrow(f(A-B)) \le \lambda^\downarrow(f(A)-f(B)).
\ee
\end{corollary}
\textit{Proof.}
Let $f=g^{-1}$, with $g$ satisfying the conditions of Proposition \ref{prop:ggc}.
Replace in (\ref{eq:ggc}) $A$ by $f(A)$ and $B$ by $f(B)$, yielding
$$
\lambda^\downarrow(g(f(A)-f(B))) \ge \lambda^\downarrow(A-B).
$$
Applying the function $f$ on both sides does not change the ordering, because of monotonicity of $f$,
and yields validity of inequality (\ref{eq:gg}).
\qed

These two results obviously imply the corresponding majorisation relations, and by Ky Fan dominance,
relations in any UI norm.
\section{Proof of Proposition \ref{prop:ggc}\label{sec:evidence}}
We want to prove inequality (\ref{eq:ggc}):
$$
\lambda^\downarrow(g(A)-g(B)) \le \lambda^\downarrow(g(A-B)),
$$
for $A,B\ge0$, $A\ge ||B||$, and concave, monotonously increasing and non-negative $g$.

W.l.o.g.\ we will assume $||B||=1$, since any other value can be absorbed in the definition of $g$.

It is immediately clear that if (\ref{eq:ggc}) holds for $g$ that in addition satisfy $g(0)=0$,
then it must also hold without that constraint, i.e.\ for functions $g(x)+c$, with $c\ge0$.
This is because the additional constant $c$ drops out in the LHS, while
$\lambda^\downarrow(g(A-B)+c)\ge \lambda^\downarrow(g(A-B))$.

Furthermore, (\ref{eq:ggc}) remains valid when replacing $g(x)$ with $ag(x)$, for $a>0$.
Thus, w.l.o.g.\ we can assume $g(0)=0$ and $g(1)=1$.
Together with concavity of $g$, this implies that, for $0\le x\le 1$, $g(x)\ge x$,
while for $x\ge 1$, the derivative $g'(x)\le 1$.

Since $0\le B\le \id$, and for $0\le x\le 1$, $g(x)\ge x$ holds, we have
$g(B)\ge B$, or $-g(B)\le -B$. By Weyl monotonicity, this implies
$\lambda^\downarrow(g(A)-g(B))\le \lambda^\downarrow(g(A)-B)$.
Thus, statement (\ref{eq:ggc}) would be implied by the stronger statement
\be
\lambda^\downarrow(g(A)-B) \le \lambda^\downarrow(g(A-B)). \label{eq:ggc3}
\ee
Now note that the argument of $g$ in the LHS, $A$, is never below 1. Thus, in principle, we could
replace $g(x)$ in the LHS by another function $h(x)$ defined as
\be
h(x) = \left\{
\begin{array}{l}
g(x),\mbox{ if }x\ge 1\\
x,\mbox{ otherwise.}
\end{array}
\right.
\ee
If we also do that in the RHS, we get a stronger statement than (\ref{eq:ggc3}). Indeed, $h(x)\le g(x)$ for $x\ge0$
and $A-B\ge0$, and therefore $h(A-B)\le g(A-B)$ holds. By Weyl monotonicity again, we see that (\ref{eq:ggc3})
is implied by
\be
\lambda^\downarrow(h(A)-B) \le \lambda^\downarrow(h(A-B)). \label{eq:ggc4}
\ee
The importance of this move is that $h(x)$ is still a monotonously increasing and concave function
(because $g'(x)\le 1$ for $x\ge1$),
but now has gradient $h'(x)\le 1$ for $x\ge0$.

Defining $C=A-B$, which is positive semi-definite,
we now have to show the inequality
$$
\lambda^\downarrow_k(h(C+B)-B) \le \lambda^\downarrow_k(h(C)) = h(\lambda^\downarrow_k(C)),
$$
for every $k$.
Fixing $k$, and introducing the shorthand $x_0=\lambda^\downarrow_k(C)$,
we can exploit concavity of $h$ to upper bound it as $h(x)\le a(x-x_0)+h(x_0)$, where
$a=h'(x_0)\le 1$.
Again by Weyl monotonicity, we find
\beas
\lambda^\downarrow_k(h(C+B)-B) &\le& \lambda^\downarrow_k(a(C+B-x_0)+h(x_0)-B) \\
&=& \lambda^\downarrow_k(aC+(a-1)B-ax_0+h(x_0)) \\
&\le& \lambda^\downarrow_k(aC)-ax_0+h(x_0) = h(x_0),
\eeas
where in the second line we could remove the term $(a-1)B$ because it is negative.
This being true for all $k$, we have proved (\ref{eq:ggc4}) and all previous statements that follow from it,
including the statement of the Theorem.
\qed
\section{On $Y$-dominated Majorisation\label{sec:domi}}
To answer Question \ref{th2}, we have to
consider the property that a convex function $f$ satisfies
\be\label{eq:star1}
\lambda(f(\Delta)) \prec_w \lambda(f(B+\Delta)-f(B))
\ee
for all PSD $B$ and $\Delta$,
which is equivalent to the statement
\be\label{eq:star2}
\lambda(f(A-B)) \prec_w \lambda(f(A)-f(B))
\ee
for all $A\ge B\ge0$.

The monotone convex angle functions $x\mapsto ax+(x-1)^+$ ($a\ge0$) already have proven their valour as a testing ground for
similar statements, in Section \ref{sec:main}.
Numerical experiments using angle functions for inequality (\ref{eq:star1}) did not directly lead
to any counterexamples, however. This temporarily increased our belief that the inequality might actually hold,
and led us to investigate, as an initial step towards a ``proof'', whether the inequality
$$
\sum_{j=1}^k \lambda_j^\downarrow(aY+B) \le \sum_{j=1}^k \lambda_j^\downarrow(aY+C)
$$
might be true for all $a\ge0$, where $B=f(Y)$ and $C=f(X+Y)-f(X)$, and $f(x)=(x-1)^+$.
The crucial observation is now that if this holds for all $a\ge0$,
then, actually, a much stronger relationship than just majorisation
must hold between $aY+B$ and $aY+C$.

To describe this phenomenon, we'll consider a somewhat broader setting.
Let $G$ and $C$ be Hermitian matrices, and let $f_1$ and $f_2$ be monotonously increasing real functions
on $\R$.
Suppose that for all $a\ge0$, the following holds:
\be\label{eq:ayb}
\sum_{j=1}^k \lambda_j^\downarrow(aA+B) \le \sum_{j=1}^k \lambda_j^\downarrow(aA+C),
\ee
with $A=f_1(G)$ and $B=f_2(G)$.

It is easily seen that if (\ref{eq:ayb}) holds for a certain value of $a$, it also holds for all smaller positive values.
Let $b$ be a scalar such that $0\le b<a$. Because both $A$ and $B$ exhibit
their eigenvalues as diagonal elements in the eigenbasis of $G$, and both in non-increasing order, we get
$$
\sum_{j=1}^k \lambda_j^\downarrow(aA+B)
=\sum_{j=1}^k \lambda_j^\downarrow(bA+B) +(a-b)\sum_{j=1}^k \lambda_j^\downarrow(A).
$$
On the other hand, for $aA+C$ we only have the subadditivity inequality
$$
\sum_{j=1}^k \lambda_j^\downarrow(aA+C)
\le\sum_{j=1}^k \lambda_j^\downarrow(bA+C) +(a-b)\sum_{j=1}^k \lambda_j^\downarrow(A).
$$
As a consequence, we obtain that, indeed,
$$
\sum_{j=1}^k \lambda_j^\downarrow(bA+B) \le \sum_{j=1}^k \lambda_j^\downarrow(bA+C)
$$
follows from (\ref{eq:ayb}).

We are therefore led to consider what happens when $a$ tends to infinity, because that value dominates all others.
Subtracting $\sum_{j=1}^k \lambda_j^\downarrow(aA)$ from both sides, and substituting $a=1/t$, we obtain
$$
\frac{1}{t}\sum_{j=1}^k (\lambda_j^\downarrow(A+tB)-\lambda_j^\downarrow(A))
\le
\frac{1}{t}\sum_{j=1}^k (\lambda_j^\downarrow(A+tC)-\lambda_j^\downarrow(A)).
$$
In the limit of $t$ going to 0, this yields a comparison between derivatives:
\be\label{eq:ddt}
\ddt\sum_{j=1}^k \lambda_j^\downarrow(A+tB)
\le
\ddt\sum_{j=1}^k \lambda_j^\downarrow(A+tC).
\ee

We will show below that the derivatives $\ddt \lambda_j^\downarrow(A+tC)$ are the diagonal
elements of $C$ in a certain basis depending on $G$ and $C$.
Let us first introduce the vector $\delta(C;A)$ whose
entries satisfy the following relation:
\be\label{eq:defdelta}
\sum_{j=1}^k \delta_j(C;A) := \ddt\sum_{j=1}^k \lambda_j^\downarrow(A+tC).
\ee
With this notation, relation (\ref{eq:ddt}) becomes
$$
\sum_{j=1}^k \delta_j(B;G) \le \sum_{j=1}^k \delta_j(C;G).
$$
That is, the entries of $\delta(B;G)$ are ``majorised'' by those of $\delta(C;G)$.
However, this is a much stronger relation than ordinary
majorisation, since the rearrangement of the entries in decreasing order is absent.

Introducing the symbol $\prec_{dw}$ for weak majorisation with missing rearrangement:
\be
a \prec_{dw} b \Longleftrightarrow \sum_{j=1}^k a_j \le \sum_{j=1}^k b_j,
\ee
relation (\ref{eq:ddt}) is expressed as
\be\label{eq:ddt2}
\delta(B;G) \prec_{dw} \delta(C;G).
\ee

To justify these notations, we now show:
\begin{proposition}\label{prop:delta}
Let $A$ and $C$ be Hermitian matrices.
With $\delta(C;A)$ defined by (\ref{eq:defdelta}),
the entries of the vector $\delta(C;A)$ are the diagonal entries of $C$
in a certain basis in which $A$ is diagonal and its diagonal entries
appear sorted in non-increasing order.
When all eigenvalues of $A$ are simple (i.e.\ have multiplicity 1), this basis is just the
eigenbasis of $A$ and does not depend on $C$.
\end{proposition}
\textit{Proof.} There are three cases to consider, according to whether $A$ is non-degenerate,
$A+tC$ has an accidental degeneracy at $t=0$, or $A+tC$ is permanently degenerate.


1. The most important case is when
all eigenvalues of $A$ are simple, i.e.\ when they have multiplicity 1.
We then show that the derivative is given by
$$
\ddt\sum_{j=1}^k \lambda_j^\downarrow(A+tC) = \trace[P_k(A)\,\, C],
$$
where $P_k(A)$ denotes the projector on the
subspace spanned by the $k$ eigenvectors of $A$ corresponding to its
$k$ largest eigenvalues.

By the simplicity of the eigenvalues of $A$, the eigenvalues of $A+tC$ are also simple for
small enough values of $t$.
This follows easily from Weyl's inequalities:
$$
\lambda_j^\downarrow(A) + \lambda_n^\downarrow(tC) \le \lambda_j^\downarrow(A+tC) \le \lambda_j^\downarrow(A) + \lambda_1^\downarrow(tC);
$$
thus if $t||C||$ is strictly less than one half the minimal difference between all pairs of eigenvalues of $A$,
the difference between all pairs of eigenvalues of $A+tC$ is bounded away from 0.
Therefore, for small enough $t$, every eigenvalue of $A+tC$ has a unique eigenvector, and as a result $P_k(A+tC)$ is well-defined as
the sum of the projectors on the eigenvectors pertaining to the $k$ largest eigenvalues.

It is well-known that the eigenvalues of $A+tC$ as functions of the real variable $t$ can be so ordered that they are
analytic functions of $t$ (see \cite{kato}, Chapter 2), and hence continuous.
This implies that the $k$-th largest eigenvalue of $A+tC$ is also a continuous function of $t$, for any $k$.

If, furthermore, an eigenvalue $\lambda(t)$ of $A+tC$ is simple in an interval of $t$,
then the projector $P(t)$ on the eigenvector $x(t)$ associated to it (with $P(t) = x(t)x(t)^*$) is also analytic, and therefore continuous
in $t$ on this interval.
We conclude that $P_k(A+tC)$ is analytic in $t$, and therefore differentiable.

By the maximality of $P_k(A)$ in the variational characterisation
$$
\sum_{j=1}^k \lambda_j^\downarrow(A) = \max_{Q_k} \trace[A\,Q_k] = \trace[A\,P_k(A)],
$$
where $Q_k$ runs over all rank-$k$ projectors, we have
$$
\ddt\trace[A\,P_k(A+tC)]=0,
$$
which implies
\beas
\lefteqn{\ddt\sum_{j=1}^k \lambda_j^\downarrow(A+tC)} \\
&=& \ddt\trace[(A+tC)\,P_k(A+tC)] \\
&=& \ddt\trace[A\,  P_k(A+tC)] + \ddt\trace[(A+tC)\,P_k(A)] \\
&=& \trace[C\,P_k(A)].
\eeas

Let $U$ be the unitary that diagonalises $A$, i.e.\ $UAU^*=\Lambda^\downarrow(A)$.
Then
$$
\trace[C\,P_k(A)] = \sum_{j=1}^k (UCU^*)_{jj},
$$
and the statement of the Proposition follows.


2. When $A$ has degenerate eigenvalues, the situation becomes somewhat more complicated, but there
are no really significant changes.
There is no longer a unique eigenbasis of $A$, so that $P_k(A)$ is not well-defined for all $k$.
We will first consider the case where $C$ is such that it removes the degeneracy of the eigenvalues
of $A$ in $A+tC$ for small enough positive $t$.
In that case $P_k(A+tC)$ will be uniquely defined for all positive $t$ less than some value $t_0$, which
is the smallest positive $t$ for which $A+tC$ has an accidental degeneracy (which is what also happens at $t=0$).

This occurs, for instance, when $C$ has simple eigenvalues.
Indeed, by analyticity of the eigenvalues of $A+tC$ in $t$, degeneracy is either accidental (for isolated values
of $t$) or permanent (for all values of $t$).
Since all eigenvalues are simple for large enough $t$, they have to remain simple for all values of $t$ except possibly for
some isolated values, such as $t=0$, in this case. Let $t_0$ be the smallest positive such value, then
$A+tC$ has simple eigenvalues for $0<t<t_0$.

We can therefore define $P_k(A)$ in a unique way as the limit $\lim_{t\to 0}P_k(A+tC)$.
This is an allowed choice because of the continuity of the eigenvalues:
$\sum_{j=0}^k \lambda^\downarrow_k(A) = \trace[\lim_{t\to 0}P_k(A+tC)\,\, A]$.
Using the same argument as in the previous case,
we obtain $\delta(C;A) := \trace[\lim_{t\to 0}P_k(A+tC)\,\, C]$.

Let $\lambda_l$ be the eigenvalues of $A$ (multiplicity not counted), and $Q_l$ the
projections onto the corresponding eigenspaces of $A$ (with $Q_l^*$ the corresponding inclusion
operators); the rank of $Q_l$ equals the multiplicity of $\lambda_l$, which we denote by $m_l$.
To obtain $\delta(C;A)$, we first construct the diagonal blocks $C_l:=Q_l C Q_l^*$ (of size $m_l$),
then take the eigenvalues $\lambda^\downarrow(C_l)$ in non-increasing order of each block, and then
concatenate the obtained sequences of eigenvalues:

$$
\delta(C;A) := (\lambda^\downarrow(C_1), \ldots, \lambda^\downarrow(C_m)).
$$
If all eigenvalues of $A$ are distinct, this reduces to the vector of diagonal elements of $C$
in the eigenbasis of $A$ that we encountered in case 1.

For example, if $\lambda^\downarrow(A)=(5,5,3,1)$, then
$\delta(C;A) = (\lambda_1^\downarrow(C_1),\lambda_2^\downarrow(C_1),C_{33},C_{44})$,
where $C_1=\twomat{C_{11}}{C_{12}}{C_{21}}{C_{22}}$ and all entries of $C$ are taken in the
eigenbasis of $A$.

Let $U$ be a unitary (which, in this case, is not unique)
that diagonalises $A$ as $UAU^*=\Lambda^\downarrow$,
and take the diagonal blocks $C_l$ of $U C U^*$, as above.
Each block can be diagonalised using a unitary $V_l$. Together with $U$ we obtain the total basis rotation
$W:=U (\bigoplus_l V_l)$. By construction, $\bigoplus_l V_l$ leaves $\Lambda$ invariant, and resolves
the ambiguity in $U$.
We obtain that $\delta(C;A)$ is the vector of diagonal entries of $C$ in the basis
obtained by applying the unitary $W$.


3. Finally, we look at the case when $A+tC$ is permanently degenerate, i.e.\
when it has degenerate eigenvalues for all values of $t$. W.l.o.g.\ we just have to look at $t$ in an interval
$[0,t_0)$, where $t_0$ is the smallest positive value for which $A+tC$ has an accidental degeneracy.
Let us denote by $\lambda_j(t)$ the eigenvalues of $A+tC$ in non-increasing order, multiplicity $m_j$ \textit{not} counted,
and by $\cP_j(t)$ the projectors on the corresponding eigenspaces.
In that case $P_k(A+tC)$ is only well-defined if there is a $j'$ such that $k=m_1+m_2+\ldots+m_{j'}$;
then we have $P_k(A+tC)=\cP_1(t)+\cP_2(t)+\ldots+\cP_{j'}(t)$.

If there is no such $j'$, let $j'$ be the largest integer such that $k> m_1+m_2+\ldots+m_{j'}=:k'$.
Thus $0<k-k'<m_{j'+1}$.
Then we have
\beas
\lefteqn{\sum_{j=1}^k \lambda^\downarrow_j(A+tC)} \\
&=& \sum_{i=1}^{j'} m_i \lambda_i(t) + (k-k')\lambda_{j'+1}(t) \\
&=& \trace[(A+tC)\,(\cP_1(t)+\ldots+\cP_{j'}(t)+\frac{k-k'}{m_{j'+1}}\cP_{j'+1}(t))] \\
&=& \trace[(A+tC)\,(\frac{k-k'}{m_{j'+1}}P_{k'+m_{j'+1}}(A+tC) + (1-\frac{k-k'}{m_{j'+1}})P_{k'}(A+tC))].
\eeas
Thus, if we define $\alpha:=(k-k')/m_{j'+1}$, $\sum_{j=1}^k \lambda^\downarrow_j(A+tC)$
interpolates linearly between $\sum_{j=1}^{k'} \lambda^\downarrow_j(A+tC)$ and
$\sum_{j=1}^{k'+m{j'+1}} \lambda^\downarrow_j(A+tC)$ with parameter $\alpha$.

Proceeding in the same way as in the two previous cases, we obtain for the derivative
$$
\ddt\sum_{j=1}^k \lambda^\downarrow_j(A+tC) = \trace[C(\alpha P_{k'+m_{j'+1}}(A) + (1-\alpha)P_{k'}(A))],
$$
where the $P_{k}(A)$ have to be replaced with the limits $\lim_{t\to 0}P_{k}(A+tC)$
if in addition there are accidental degeneracies at $t=0$.
Let us consider the entries of $C$ again as before, in an eigenbasis of $A$ in which the eigenvalues of $A$
appear on the diagonal, in non-increasing order.
We get
\be\label{eq:interpol}
\delta(C;A)_k = (1-\alpha) \sum_{i=1}^{k'} C_{ii} +\alpha\sum_{i=1}^{k'+m_{j'+1}} C_{ii}.
\ee
Because of the permanent degeneracy, an eigenbasis
is determined up to ``local'' rotations within the various eigenspaces.
We consider a partitioning of $C$ in such an eigenbasis corresponding to these eigenspaces. That is,
in $C$ we can single out diagonal blocks, each of which corresponds to the eigenspace of eigenvalue $\lambda_j$.
We can use our freedom to choose the local rotations to make all diagonal elements of $C$ equal within
each diagonal block. This allows us to get rid of the interpolation in (\ref{eq:interpol}), and we finally obtain
that, again,
$$
\delta(C;A)_k = \sum_{i=1}^{k} C_{ii},
$$
with the entries of $C$ taken in the eigenbasis that we have just chosen.
\qed

The upshot of this Proposition is that there exists a unitary $U$ such that
$UAU^*=\Lambda^\downarrow(A)$ and $\delta(C;A)=\diag(UCU^*)$. In the generic case
that all $\lambda_i(A)$ are distinct, $U$ is unique and does not depend on $C$.

A number of easy consequences follow immediately from this Proposition:
\begin{corollary}
Let $G$ and $C$ be Hermitian matrices, $f$ be any monotonously increasing real function on $\R$, and $g$ any
strictly increasing real function on $\R$, then
\begin{enumerate}
\item[(i)] $\delta(f(G);G) = f(\lambda^\downarrow(G))$.
\item[(ii)] $\delta(C;G)$ obeys Schur's majorisation Theorem: $\delta(C;G)\prec \lambda^\downarrow(C)$.
\item[(iii)] $\delta(C;G) + a\lambda^\downarrow(f(G)) = \delta(C+af(G);G), \forall a\ge0$.
\item[(iv)] $\delta(C;f(A))=\delta(C;A)$.
\end{enumerate}
\end{corollary}

Along with the previously demonstrated equivalence of (\ref{eq:ayb})
with (\ref{eq:ddt2}), the Corollary immediately leads to the following Proposition:
\begin{proposition}\label{prop:4b}
For Hermitian $G,C$, monotonously increasing real functions $f_1,f_2$ on $\R$,
and $A=f_1(G)$, $B=f_2(G)$,
the following are equivalent:
\bea
\lambda(aA+B) &\prec_w& \lambda(aA+C), \quad\forall a\ge0 \label{eq:pp1}\\
\delta(B;G) &\prec_{dw}& \delta(C;G) \label{eq:pp2} \\
\delta(aA+B;G) &\prec_{dw}& \delta(aA+C;G), \quad\forall a\ge0.\label{eq:pp3}
\eea
\end{proposition}
\textit{Proof.}

(\ref{eq:pp1}) implies (\ref{eq:pp2}): This is just Proposition \ref{prop:delta}.

(\ref{eq:pp2}) implies (\ref{eq:pp3}): Add $a\lambda^\downarrow(A)$ to both sides and invoke statement (iii) of the Corollary.

(\ref{eq:pp3}) implies (\ref{eq:pp1}):
By statement (i) of the Corollary, the LHS of (\ref{eq:pp3}) is equal to $\lambda^\downarrow(a A+B)$,
while, by statement (ii) of the Corollary, its RHS is majorised by $\lambda(aA + C)$.
\qed
\section{Counterexample to Question \ref{th2}\label{sec:counter}}
If the answer to Question \ref{th2} is to be affirmative, it should at least
hold for all angle functions $f(x)=ax+b(x-x_0)^+$.
By Proposition \ref{prop:4b} this is equivalent to the statement
$$
\delta((Y-\id)^+;Y) \prec_{dw} \delta((X+Y-\id)^+ - (X-\id)^+;Y).
$$
Consider the $3\times 3$ matrices
$$
X = \left(
\begin{array}{rrr}
         0.35614  &       -0.053243   &      0.10116\\
         -0.053243  &       0.87456   &     0.40559\\
         0.10116  &      0.40559   &      0.82474
\end{array}
\right)
$$
and
$$
Y = \left(
\begin{array}{rrr}
         0.53642  &                       0   &                      0\\
                         0  &       0.42018   &                      0\\
                         0  &                       0   &      0.094866
\end{array}
\right).
$$
The eigenbasis of $Y$ is therefore the standard basis.
Then $\delta((Y-\id)^+;Y)=(0,0,0)$
and
$$
(X+Y-\id)^+ - (X-\id)^+ = \left(
\begin{array}{rrr}
-0.00018194  & 0.00052449 &  -0.0016345\\
   0.00052449 &      0.2573 &     0.12368\\
   -0.0016345  &    0.12368 &        0.04
\end{array}
\right)
$$
so that $\delta((X+Y-\id)^+ - (X-\id)^+;Y)=(-0.00018194, 0.2573, 0.04)$.
The first entry is negative, violating the $\prec_{dw}$ relation,
and thereby answering Question \ref{th2} in the negative.
\qed
\section{Further Applications of $Y$-dominated majorisation}
One issue we had to address during our attempts at giving a positive answer to Question \ref{th2}
dealt with the possibility of reducing the question for convex functions to convex angle functions.
One way of doing so would have been possible if the set of (monotonously increasing and convex) functions
satisfying (\ref{eq:star1}) were closed under addition.
While we were unable to prove this particular statement (which is most likely false, anyway),
Proposition \ref{prop:4b} enables us to prove the corresponding statement for the relation
\be\label{eq:star3}
\delta(f(Y);Y) \prec_{dw} \delta(f(X+Y)-f(X);Y).
\ee
\begin{proposition}
Let all the eigenvalues of $Y$ be distinct.
Let $f$ and $g$ be functions from $\R$ to $\R$ satisfying (\ref{eq:star3}).
Then $f+g$ also satisfies (\ref{eq:star3}).
\end{proposition}
\textit{Proof.}
By the assumption on the eigenvalues of $Y$, $\delta(A;Y)$ equals $\diag(A)$ in a basis only depending on $Y$
and is therefore a linear function of $A$.
We can therefore add up the inequalities (\ref{eq:star3}) for $f$ and $g$ and obtain
the corresponding inequalities for $f+g$.
\qed

\bigskip

A second application of Proposition \ref{prop:4b} is a strengthening of the following Proposition,
which we also obtained in the course of our attempts at positively answering Question \ref{th2}.
\begin{proposition}\label{prop:g}
For $X,Y\ge0$ and $g_a(x)=ax+\frac{x^2}{x+1}$, with $a\ge0$, the following majorisation statement holds:
$$
\lambda(g_a(Y)) \prec_w \lambda(g_a(X+Y)-g_a(X)).
$$
\end{proposition}
\textit{Proof.}
From the proof of Lemma X.1.4 in \cite{bhatia}, we have, for $X,Y\ge0$,
$$
\lambda_j^\downarrow( (X+\id)^{-1} - (X+Y+\id)^{-1} ) \le
\lambda_j^\downarrow( \id - (Y+\id)^{-1} ).
$$
Defining the function $f(x)=\frac{x}{x+1}=1-(x+1)^{-1}$,
this turns into:
$$
\lambda_j^\downarrow(f(X+Y)-f(X)) \le \lambda_j^\downarrow(f(Y)).
$$
This implies the majorisation statement
\be\label{eq:maj1}
\sum_{j=1}^k  \lambda_j^\downarrow(f(X+Y)-f(X)) \le \sum_{j=1}^k \lambda_j^\downarrow(f(Y)).
\ee

We want to prove a somewhat similar statement for the function $g_a(x)$.
Since both $f$ and $g_a$ are monotonously increasing over $\R^+$, and noting that $g_a(x)= (a+1)x-f(x)$, we have
\beas
\lambda_j^\downarrow(g_a(Y)) &=& g_a(\lambda_j^\downarrow(Y))
= (a+1)\lambda_j^\downarrow(Y) - f(\lambda_j^\downarrow(Y)) \\
\lambda_j^\downarrow(f(Y)) &=& f(\lambda_j^\downarrow(Y)),
\eeas
so that
$$
\lambda_j^\downarrow(g_a(Y)) = (a+1)\lambda_j^\downarrow(Y) - \lambda_j^\downarrow(f(Y)).
$$
This implies in particular
\beas
\sum_{j=1}^k \lambda_j^\downarrow(g_a(Y)) &=&
(a+1)\sum_{j=1}^k \lambda_j^\downarrow(Y) - \sum_{j=1}^k \lambda_j^\downarrow(f(Y)) \\
&\le& (a+1)\sum_{j=1}^k \lambda_j^\downarrow(Y) - \sum_{j=1}^k  \lambda_j^\downarrow(f(X+Y)-f(X)),
\eeas
where we have inserted (\ref{eq:maj1}).
Exploiting the well-known relation (\cite{bhatia}, Th.\ III.4.1)
$$
\sum_{j=1}^k \lambda_j^\downarrow(A+B)
\le \sum_{j=1}^k \lambda_j^\downarrow(A) + \sum_{j=1}^k \lambda_j^\downarrow(B),
$$
for $A=(a+1)Y-f(X+Y)+f(X)$ and $B=f(X+Y)-f(X)$ then yields
\beas
\sum_{j=1}^k \lambda_j^\downarrow(g_a(Y)) &\le& \sum_{j=1}^k \lambda_j^\downarrow((a+1)Y-f(X+Y)+f(X)) \\
&=& \sum_{j=1}^k \lambda_j^\downarrow(g_a(X+Y)-g_a(X)).
\eeas
\qed

Proposition \ref{prop:4b}, with $A=G=Y$,
$B=f(Y)$, $C=f(X+Y)-f(X)$, where $f(x)=x^2/(x+1)$, then yields
the following strengthening of Proposition \ref{prop:g}:
\begin{proposition}\label{prop:4}
For $X,Y\ge0$, and $g_a(x)=ax+\frac{x^2}{x+1}$, with $a\ge0$,
$$
\delta(g_a(Y);Y) \prec_{dw} \delta(g_a(X+Y)-g_a(X);Y).
$$
\end{proposition}
Here we noted that $g_a(X+Y)-g_a(X) = aY+f(X+Y)-f(X)$.

\bigskip

To end this Section, we present a third application of Proposition \ref{prop:4b},
namely to the results of Kosem and Bourin-Uchiyama.
Consider first inequality (\ref{eq:AZ1}), which holds for all non-negative concave functions $f(x)$.
In particular, it holds for all functions $f=ax+f_0(x)$, where $f_0$ is non-negative concave, and $a\ge0$.
Inserting this in the eigenvalue-majorisation form of inequality (\ref{eq:AZ1}), we get the
$(A+B)$-dominated majorisation relation
$$
\lambda(a(A+B)+f_0(A+B))\prec_w \lambda(a(A+B)+f_0(A)+f_0(B)),
$$
for $A,B\ge0$.
Proposition \ref{prop:4b} then immediately yields the stronger form
\be\label{eq:bourins}
\delta(f(A+B);A+B)\prec_{dw} \delta(f(A)+f(B);A+B),
\ee
for all non-negative concave functions $f$.
The strengthening of inequality (\ref{eq:AZ2}) is performed in a completely identical way and yields
the reversed inequality of (\ref{eq:bourins}) for non-negative convex functions $g$ such that $g(0)=0$.
\begin{ack}
JSA thanks Professor Moin Uddin, Director of his institute for encouragement and
supporting his visit to attend the conference at Nova Southeastern University, Fort
Lauderdale, Florida, USA, which lead to his introduction to
Koenraad  M.R.~Audenaert and the completion of this work.

KA thanks the Institute for Mathematical Sciences, Imperial College London, for support.
His work is part of the QIP-IRC (www.qipirc.org) supported by EPSRC (GR/S82176/0).
\end{ack}

\end{document}